\documentclass[a4paper,11pt,oneside,notitlepage]{article}
\usepackage{times}
\usepackage{lmodern}
\usepackage[T1]{fontenc}
\usepackage[utf8]{inputenc}
\usepackage[english]{babel}
\usepackage{amssymb}
\usepackage{amsmath}
\usepackage{amsthm}
\usepackage{amsfonts}
\usepackage{latexsym}
\usepackage[pdftex]{color, graphicx}
\usepackage{makeidx}
\usepackage{sidecap}
\usepackage{verbatim}
\usepackage{geometry}
\usepackage{subfig}
\usepackage{tikz,float}
\usetikzlibrary{shapes,arrows,shadows, svg.path}
\usepackage{caption}
\definecolor{mygray}{RGB}{201,201,201}

\newtheorem{thm}{Theorem}[section]

\newtheorem{cor}[thm]{Corollary}

\def\laplace{\Delta}
\def\utheta{U_{\Theta}^{\varepsilon}}

\def\vtheta{V_{\Theta}^{\varepsilon}}

\def\normalder{\partial_{\nu}}
\def\ballzero{\Omega_0}
\def\ballone{\Omega_1}

\def\eigenv{\lambda^{\varepsilon}}
\def\grad{\nabla_x}
\def\eiGelf{\Lambda^{\varepsilon}}

\def\perceltheta{\omega_{\Theta}}
\def\percelupsilon{\omega_{\Upsilon}}
\def\unitsquare{Q}
\def\ucomptheta{U^\varepsilon_{\unitsquare\setminus\Theta}}
\def\vcomptheta{V^\varepsilon_{\unitsquare\setminus\Theta}}

\def\uone{u_1^{\varepsilon}}
\def\uzero{u_0^{\varepsilon}}
\def\espec{\sigma_{\text{e}}}
\def\ball{\Theta}
\def\complball{\unitsquare\setminus\ball}
\def\perfunc{\mathcal{H}_\eta(\perceltheta)}

\numberwithin{equation}{section}

\title{Gaps in the spectrum of two-dimensional square packing of stiff disks}
\author{L. D'Elia \footnote{Dipartimento di Scienze Matematiche, Politecnico di Torino
 Corso Duca degli Abruzzi 24,
10129 Torino, Italy {\tt e-mail: lorenza.delia@polito.it}}
 \ and S.A. Nazarov \footnote{St. Petersburg State University, Universitetskaya nab., 7-9, St. Petersburg, 199034, Russia\hspace{0.2cm} \ and \hspace{0.2cm}
 Institute of Problems Mechanical Engineering RAS, V.O., Bolshoj pr., 61, St. Petersburg, 199178, Russia {\tt e-mail: srgnazarov@yahoo.co.uk}}}

\date{}

\begin{document}
\maketitle

\abstract{\noindent In this paper we investigate via an asymptotic method the opening  of gaps in the spectrum of a stiff problem for the Laplace operator $-\Delta$ in $\mathbb{R}^2$ perforated by contiguous circular holes. The density and the stiffness constants are of order $\varepsilon^{-2m}$ and $\varepsilon^{-1}$ in the holes with $m\in (0,1/2)$. We provide an explicit expression of the leading terms of the eigenvalues and the corresponding eigenfunctions which are related to the Bessel functions of the first kind.  }


\section{Introduction}
\subsection{Formulation of problem}
\begin{figure}[t]
\centering
\begin{minipage}{0.3\textwidth}
\centering
\includegraphics[scale=0.22]{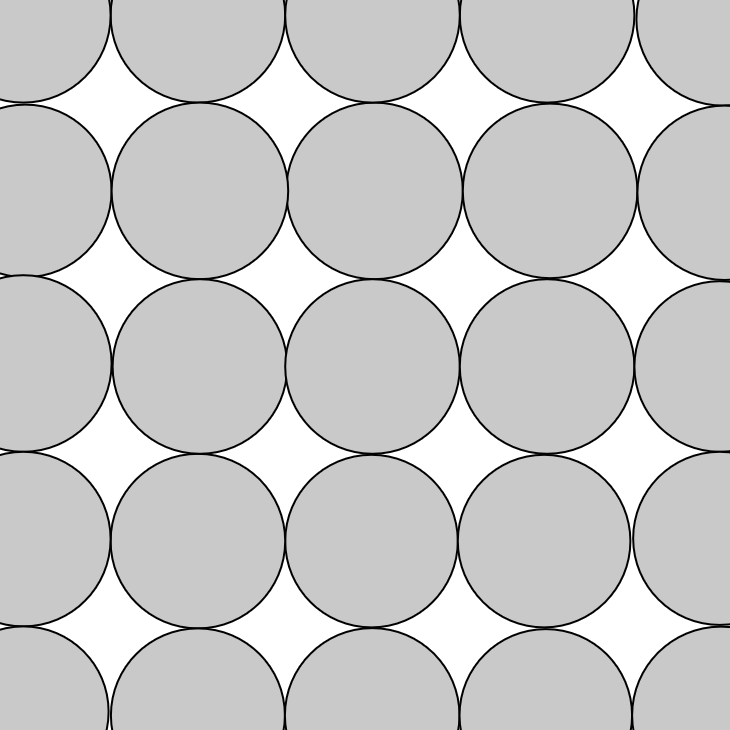}
\caption*{(a)}
\end{minipage} %
\begin{minipage}{0.3\textwidth}
\centering
\includegraphics[scale=0.55]{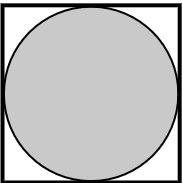}
\caption*{(b)}
\end{minipage} %
\hspace{0.5cm}
\begin{minipage}{0.3\textwidth}
\centering
\includegraphics[scale=0.55]{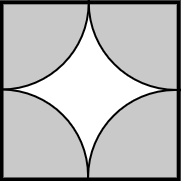}
\caption*{(c)}
\end{minipage}
 \caption{\small Figure (a) shows $\Omega_0\cup\Omega_1$. Two (possible) choices $\perceltheta$ and $\percelupsilon$ of the periodicity cell are drawn in Figures (b) and (c) respectively} 
 \label{fig:periodicells}
\end{figure}
\noindent Let $\ballzero$ be the plane $\mathbb{R}^2$ perforated by contiguous circular holes 
      \begin{equation*}
      B_{1/2} (\alpha) :=\{ x=(x_1,x_2)\hspace{0.02cm}:\hspace{0.02cm} (x_1-\alpha_1, x_2-\alpha_2)\in B_{1/2}\},
      \end{equation*}
where $\alpha=(\alpha_1,\alpha_2)\in\mathbb{Z}^2$ is a multi-index, $\mathbb{Z}=\{0,\pm 1, \pm2,\dots\}$ and $ B_{1/2} :=\{x\hspace{0.02cm}:\hspace{0.02cm} |x|< 1/2\}$. More precisely, 
     \begin{equation*}
     \ballzero := \mathbb{R}^2\setminus \bigcup _{\alpha\in\mathbb{Z}^2}\overline{B_{1/2} (\alpha)}. 
     \end{equation*}
We set 
     \begin{equation*}
     \ballone :=\bigcup _{\alpha\in\mathbb{Z}^2}B_{1/2} (\alpha)\quad\mbox{and}\quad \partial\ballone := \bigcup  _{\alpha\in\mathbb{Z}^2}\partial B_{1/2} (\alpha).
     \end{equation*}
We consider the stiff spectral problem in the inhomogeneous plane (see Figure \ref{fig:periodicells}(a))
     \begin{align}
            -\laplace u^{\varepsilon}_1(x) &= \eigenv u^{\varepsilon}_1(x), \hspace{4cm} x\in\ballone,\label{utheta} \\
            -\varepsilon^{-1}\laplace u^{\varepsilon}_0(x) &= \varepsilon^{-2m}\eigenv u^{\varepsilon}_0(x), \hspace{3.2cm} x\in\ballzero, \notag\\
            u^{\varepsilon}_1(x) = u^{\varepsilon}_0(x), &\qquad \varepsilon^{-1}\normalder u^{\varepsilon}_0(x) = \normalder u^{\varepsilon}_1(x),\hspace{1.5cm} x\in \partial\ballone,\label{sbc}
           \end{align}
Here, $\laplace$ is the Laplace operator, $\eigenv$  the spectral parameter, $\nu$  the outward unit normal vector to $\partial\Omega_1$, $\normalder = \nu\cdot\nabla$  the normal derivative, $\nabla$  the gradient and $m\in (0, 1/2)$  stands for a fixed exponent. 
For any $\varepsilon>0$, the variational setting of the problem \eqref{utheta}-\eqref{sbc} reads as
    \begin{equation}
    \label{varform}
    (\nabla \uone, \nabla\varphi_1)_{\ballone} +\varepsilon^{-1}(\nabla \uzero, \nabla\varphi_0)_{\ballzero} = \eigenv \left((\uone,\varphi_1)_{\ballone} + \varepsilon^{-2m}(\uzero,\varphi_0)_{\ballzero}\right), 
    \end{equation}
for $\varphi\in H^1( \mathbb{R}^2)$, 
where $(\cdot, \cdot)_{\Omega_j}$ stands for the natural inner product in $L^2(\Omega_j)$, for $j=0,1$.  
To the problem \eqref{varform} we  assign a positive and self-adjoint operator $A_{\varepsilon}$ in the Hilbert space $L^2(\mathbb{R}^2)$ with the  domain  $\mathcal{D}(A_{\varepsilon})\subset H^1(\mathbb{R}^2)$ (see \cite[Ch. 10]{BS87}). Furthermore, owing to general results in \cite{N08Neu}, see also \cite[Section 8]{CDN}, we have 
      \begin{equation*}
      {\cal D}(A_\varepsilon) := \{u^\varepsilon \hspace{0.02cm}:\hspace{0.02cm} u^\varepsilon_j\in H^2(\Omega_j),\hspace{0.05cm} \mbox{for } j=1,2, \hspace{0.05cm}\mbox{and }\eqref{sbc} \mbox{ holds} \},
      \end{equation*}
that is, tips of cusps do not bring serious singularities to eigenfucntions of problem \eqref{utheta}-\eqref{sbc}.
The spectrum $\sigma$ of $A_\varepsilon$ is contained in the positive semi-axis $\overline{\mathbb{R}}_{+}:= [0,\infty)$ and  since the embedding $H^1(\mathbb{R}^2)\subset L^2(\mathbb{R}^2)$ is not compact, the essential spectrum $\espec^\varepsilon$ does not consist of the single point $\lambda=0$ (see \cite[Theorem 10.1.5]{BS87}). Such an essential spectrum $\espec^\varepsilon$ has a band-gap structure (see, \textit{e.g.}, \cite{K82, RS78, S85}), \textit{i.e.} it is  represented as the countable union
     \begin{equation}
     \label{essespe}
     \espec^\varepsilon=\bigcup_{n=1}^{\infty} {\cal B}^\varepsilon_n,
     \end{equation}
of the compact and connected spectral bands 
          \begin{equation}
          \label{segment}
          {\cal B}^\varepsilon_n:=\{\lambda^\varepsilon_n = \Lambda^\varepsilon_n(\eta)\hspace{0.02cm}|\hspace{0.05cm} \eta\in [-\pi,\pi)^2\}.
          \end{equation}
The bands ${\cal B}^\varepsilon_n$ involve entries of monotone increasing unbounded positive sequence
        \begin{equation}
        \label{spectrum}
        0\leq\Lambda_1^\varepsilon(\eta) \leq \Lambda_2^\varepsilon(\eta)\leq\dots\leq\Lambda^\varepsilon_n(\eta)\leq\dots\to\infty
        \end{equation}
of eigenvalues of the auxiliary spectral problem on the periodicity cell 
$\perceltheta:=\ball \cup (\unitsquare\setminus\ball)$
    \begin{align}
    -(\grad +i\eta)^2\utheta(x,\eta) &= \eiGelf(\eta) \utheta(x,\eta), \hspace{1.7cm} x\in\ball,\label{uthetaG} \\
    -\varepsilon^{-1}(\grad +i\eta)^2 \ucomptheta(x,\eta) &= \varepsilon^{-2m}\eiGelf(\eta)\ucomptheta(x,\eta), \hspace{0.5cm} x\in\complball, \label{uupsiG}\\
    \utheta(x, \eta) &= \ucomptheta(x,\eta), \hspace{2.2cm} x\in \Gamma, \notag\\
    \varepsilon^{-1}\nu\cdot(\grad+i\eta) \ucomptheta(x,\eta) &= \nu\cdot(\grad+i\eta) \utheta(x,\eta), \hspace{0.7cm} x\in \Gamma,\label{sbcG}
    \end{align}
along with the periodicity conditions
    \begin{align}
     \ucomptheta(\tfrac{1}{2}, x_2,\eta) =  \ucomptheta(-\tfrac{1}{2}, x_2,\eta),&\hspace{0.7cm}  \ucomptheta( x_1,\tfrac{1}{2},\eta) =  \ucomptheta(x_1,-\tfrac{1}{2},\eta),\label{quasipercond1}\\
    \frac{\partial  }{\partial x_1}\ucomptheta(\tfrac{1}{2}, x_2,\eta) =\frac{\partial  }{\partial x_1}\ucomptheta (-\tfrac{1}{2}, x_2,\eta),&\hspace{0.7cm}  \frac{\partial  }{\partial x_2}\ucomptheta( x_1,\tfrac{1}{2},\eta) =\frac{\partial  }{\partial x_2}\ucomptheta( x_1,-\tfrac{1}{2},\eta) \label{quasipercond2},
    \end{align}
where $\unitsquare := (-1/2, 1/2)^2$ is the unit square in $\mathbb{R}^2$, $\ball:= B_{1/2}$ is the disk inside $\unitsquare$  and $\Gamma := \partial\ball$ (see Figure \ref{fig:periodicells}(b)). Here, $\grad$ denotes the gradient with respect to the variable $x$ and $i$ stands for the imaginary unity. The multiplicities of the eigenvalues in \eqref{spectrum} is taken into account and the functions $U^\varepsilon_\Theta$ and $U^\varepsilon_{Q\setminus\Theta}$ are  the Gelfand images of  $u^\varepsilon_1$ and $u^\varepsilon_0$ respectively, where the Gelfand transform (see \cite{G50}), also known as the Floquet-Bloch transform (see  \cite{K82, K93, S85}), is defined by 
    \begin{equation}
    \label{Gelfandtran}
    u(x)\mapsto U(x,\eta) := \frac{1}{2\pi}\sum_{k\in\mathbb{Z}^2}e^{-i\eta\cdot(x+k)}u(x+k),
    \end{equation}
with $\eta\in [-\pi, \pi)^2$ being the Floquet parameter. Note that the variable $x$ on the left-hand side of \eqref{Gelfandtran} belongs to  $\mathbb{R}^2$ but on the right-hand side $x$ lives in the periodicity cell $\perceltheta$.  If $u\in C^\infty_{\rm c}(\mathbb{R}^2)$, $(\partial_{x_j}+i\eta)u(x)$ is the Gelfand image  of the partial derivative $\partial_{x_j}$, for $j=1,2$, so that, applying the transform \eqref{Gelfandtran} to problem \eqref{utheta}-\eqref{sbc}, we obtain the parameter $\eta$-dependent problem \eqref{uthetaG}-\eqref{quasipercond2} in the cell $\perceltheta$. 
\par For any $\eta\in [-\pi,\pi)^2$, the problem \eqref{uthetaG}-\eqref{quasipercond2} is associated with a positive and self-adjoint operator $A^\varepsilon(\eta)$. The advantage of dealing with the family of operators $A^\varepsilon(\eta)$ is that such a family has discrete spectrum given by  \eqref{spectrum}, since the embedding $\perfunc\subset L^2(\perceltheta)$ is compact, where $\perfunc$ denotes the space of functions in $H^1(\perceltheta)$ satisfying the periodicity condition \eqref{quasipercond1}. It is known (see \textit{e.g.} \cite[Chapter 6]{K95} and \cite[Chapter 9]{K93}) that the functions $$\eta\in [-\pi,\pi)^2 \mapsto \eiGelf_n(\eta)$$ are continuous  and $2\pi$-periodic, so that, the spectral bands \eqref{segment} are  compact real intervals. 
\par  The aim of the present paper is to study the band-gap structure of the spectrum \eqref{essespe} of the problem \eqref{utheta}-\eqref{sbc} and to discuss the opening of the spectral gaps ${\cal G}^\varepsilon$ through an asymptotic method. Recall that the gaps ${\cal G}^\varepsilon_n$ are open intervals free of the essential spectrum $\espec^\varepsilon$  with endpoints in the $\espec^\varepsilon$. These gaps occur when the bands ${\cal B}^\varepsilon_n$ do not overlap and touch each other.  To study the existence of ${\cal G}^\varepsilon_n$, the results obtained in \cite{CDN} are crucial.  Indeed, the formal ans\"{a}tze of the eigenpairs $(\Lambda^\varepsilon, \{U^\varepsilon_\Theta, U^\varepsilon_{Q\setminus\Theta}\})$ of the model problem \eqref{uthetaG}-\eqref{sbcG} along with \eqref{quasipercond1}-\eqref{quasipercond2} are suggested by the ones performed in \cite{CDN}. We point out that in \cite{CDN} the authors have discussed the expansions of the eigenpairs for any value of the parameter $m\in\mathbb{R}$. In our analysis, we focus only on the case $m\in (0, 1/2)$, since this case seems more interesting due to difficulties arising from the dependence of $\eta\in [-\pi, \pi)^2$. At the same time, we demonstrate all necessary tools to adapt the analysis to other values of $m$. 
\par In the present context, the computation of the terms appearing in the asymptotic expansions is trickier than that done in \cite{CDN}. The main issue is related to the geometry of the periodicity cell $\perceltheta$. More specifically, the non-connectedness of $Q\setminus\Theta$ in the periodicity cell $\perceltheta$ makes difficult the explicit computation of  the leading and  first-order correction terms of the asymptotic expansion of $U^\varepsilon_{Q\setminus\Theta}$. To overcome this obstacle, we exploit the geometry of the inhomogeneous plane so that  we choose another version of the periodicity cell where $Q\setminus\Theta$ turns into a connected domain. Such a periodicity cell $\percelupsilon$ is drawn in Figure \ref{fig:periodicells}(c) and is given by $\Upsilon\cup(\unitsquare\setminus\Upsilon)$, where $\Upsilon$ is defined by
    \begin{equation*}
    \Upsilon := \{x\in Q\hspace{0.02cm}:\hspace{0.02cm} |x -P_{j\pm}|>1,\hspace{0.3cm} j=1,2   \},
    \end{equation*} 
where $P_{j\pm}$ are the vertices of unit square $\unitsquare$, {\it i.e.} $P_{1\pm} := (+ 1/2,\pm 1/2)$, and $P_{2\pm} := (- 1/2,\pm 1/2)$. 
In other words, $\percelupsilon$ is obtained by eliminating a quarter of the unit disks centered at the vertices $P_{j\pm}$, for $j=1,2$, and radius $1/2$ from the unit square $Q$. Therefore, due to the periodicity conditions \eqref{quasipercond1}, the disconnected set $Q\setminus\Upsilon$ in $\perceltheta$ turns into the connected domain $\Upsilon$ in $\percelupsilon$ and hence boundary value problems in $\percelupsilon$ are solved in the classical Sobolev spaces (see Sections 2.2 and 2.3).
\par As far as the justification procedure is concerned, we follow the same lines of the proof of \cite[Theorem 3.1]{CDN}. We enlighten that in \cite{CDN} the authors handle a different geometry of domain and a problem where the Laplace operator combined with a \textquotedblleft classical\textquotedblright\, normal derivative appears, so that, at first glance, the problem in \cite{CDN} seems different from the model problem \eqref{uthetaG}-\eqref{sbcG} studied in our context. However, we can obtain a similar problem using some tricks. Indeed, in light of the geometry of the inhomogeneous plane, we may mainly choose the periodicity cell in two ways depending on the position of the cusp points. Such cusps  may lie on the boundary of the periodicity cell where the periodicity conditions are imposed, such as $\perceltheta$ and $\percelupsilon$, or the cusp points are in the interior of the unit square $Q$, such as $\omega$ (see Figure \ref{figure:cellc}). The latter choice allows us to recover the same geometry as in \cite{CDN}.  More specifically, we translate the unit square $Q$ of the vector $(1/4, 3/4)$ and we denote by $Q'$ the translated square, {\it i.e.} $Q':= (-1/4, 3/4)\times (1/4, 5/4)$. Hence, the periodicity cell $\omega$ is given by  $\Xi\cup (\omega\setminus\Xi)$, where $\Xi:= Q'\cap\Omega_1$. 
\par In order to obtain a problem where the Laplace operator appears, we introduce an equivalent version of Gelfand's transform  given by
    \begin{equation}
    \label{newGelftransint}
    u(x)\to G(x,\eta) := {1\over 2\pi}\sum_{k\in\mathbb{Z}^2}e^{-i\eta\cdot k} u(x+k).
    \end{equation}
Applying \eqref{newGelftransint} to problem \eqref{utheta}-\eqref{sbc}, the model problem in $\omega$ turns into
\begin{align}
    -\laplace G_\Xi^\varepsilon(x,\eta) &=\Lambda^\varepsilon(\eta) G_\Xi^\varepsilon(x,\eta), \hspace{1.5cm} x\in \Xi,\label{modlpb}\\
     -\varepsilon^{-1}\laplace G_{\omega\setminus\Xi}^\varepsilon(x,\eta) &=\varepsilon^{-2m}\Lambda^\varepsilon(\eta) G_{\omega\setminus\Xi}^\varepsilon(x,\eta), \hspace{0.3cm} x\in \omega\setminus\Xi,\notag\\
    G_\Xi^\varepsilon(x,\eta) &=G_{\omega\setminus\Xi}^\varepsilon(x,\eta), \hspace{1.9cm}x\in\Gamma_\Xi, \notag\\
    \varepsilon^{-1}\partial_\nu  G_{\omega\setminus\Xi}^\varepsilon(x,\eta)&=\partial_\nu  G_{\Xi}^\varepsilon(x,\eta), \hspace{1.9cm}x\in\Gamma_\Xi,\label{modpb1}
    \end{align}
together with the quasi-periodicity conditions
  \begin{align}
     G_{j}^\varepsilon(\tfrac{3}{4}, x_2,\eta) = e^{i\eta_1} G_{j}^\varepsilon(-\tfrac{1}{4}, x_2,\eta),&\hspace{0.7cm}  G_{j}^\varepsilon( x_1,\tfrac{5}{4},\eta) = e^{i\eta_2} G_{j}^\varepsilon(x_1,\tfrac{1}{4},\eta),\label{qptrue1}\\
    \frac{\partial  }{\partial x_1}G_{j}^\varepsilon(\tfrac{3}{4}, x_2,\eta) &=e^{i\eta_1}\frac{\partial  }{\partial x_1}G_{j}^\varepsilon (-\tfrac{1}{4}, x_2,\eta),\notag\\ 
    \frac{\partial  }{\partial x_2}G_{j}^\varepsilon( x_1,\tfrac{5}{4},\eta) &=e^{i\eta_2}\frac{\partial  }{\partial x_2}G_{j}^\varepsilon( x_1,\tfrac{1}{4},\eta),\label{qptrue3}
    \end{align}
for $j=\Xi, \omega\setminus\Xi$, where $G_\Xi^\varepsilon$ and $G_{\omega\setminus\Xi}^\varepsilon$ are the image through the Gelfand transform \eqref{newGelftransint} of $u^\varepsilon_1$ and $u^\varepsilon_0$ respectively. Here $\Gamma_\Xi$ is the boundary of $\Xi$. For any $\eta\in [-\pi, \pi)^2$, we assign to problem  \eqref{modlpb}-\eqref{qptrue3} a positive and self-adjoint operator $A_{\rm qp}^\varepsilon(\eta)$. Since the embedding $H^1_{\rm qp}\subset L^2(\omega)$ is compact, where $H^1_{\rm qp}$ is the subspace of the Sobolev space $H^1(\omega)$ satisfying the quasi-periodicity conditions \eqref{qptrue1}, the spectrum of $A^\varepsilon_{\rm qp}(\eta)$ is given by the sequence \eqref{spectrum}. Then,   
the model problem \eqref{modlpb}-\eqref{qptrue3} in the periodicity cell $\omega$ enables us to repeat the same arguments of \cite[Theorem 3.1]{CDN} whose claim is stated here for the convenience of the reader (see Theorem \ref{thmmain}).  We point out that the Floquet parameter $\eta$ does not bring a trouble because in the last setting   all attributes of the model problem, including constants in {\it a priori} estimated, depend on $\eta$ continuously (see, {\it e.g.} \cite{K95}).
\par The main feature of the present paper is the explicit computation of the leading terms $\Lambda^0$ and $U^0_\Theta$ of ans\"{a}tze of $\Lambda^\varepsilon$ and $U^\varepsilon_\Theta$.  This is a consequence of the geometry of the inhomogeneous plane and of the periodicity cell $\perceltheta$. Indeed, the limit problem in the disk $\Theta$ turns out to be the Helmholtz equation combined with the homogeneous Dirichlet condition on $\Gamma$, so that  $U^0_{\Theta}$ and $\Lambda^0$ are  expressed through the Bessel functions $J_n$ of the first kind  and their zeros $j_{n,k}$ for $n\in\mathbb{N}\cup \{0\}$ and $ k \in\mathbb{N}$, where $\mathbb{N}$ is the set of positive integer number $\{1,2,3\dots\}$. Therefore,  the spectrum of the limit problem in the disk $\Theta$ is given by the monotone increasing, unbounded and positive sequence  of the eigenvalues
\begin{equation}
\label{spectrumlimprob}
       0<\Lambda^0_{0,1}<\Lambda^0_{1,1c}=\Lambda^0_{1,1s}  <\Lambda^0_{2,1c}=\Lambda^0_{2,1s}  <\Lambda^0_{0,2}<\Lambda^0_{3,1c}=\Lambda^0_{3,1s} < \Lambda^0_{1,2c}=\Lambda^0_{1,2s}<\dots,
       \end{equation}
which implies that  the sequence \eqref{spectrum} becomes
\begin{align}
       \Lambda^\varepsilon_{0,1}(\eta)<\Lambda^\varepsilon_{1,1c}(\eta)&\leq\Lambda^\varepsilon_{1,1s} (\eta) <\Lambda^\varepsilon_{2,1c}(\eta)\leq\Lambda^\varepsilon_{2,1s} (\eta) <\Lambda^\varepsilon_{0,2}(\eta)<\Lambda^\varepsilon_{3,1c}(\eta)\notag\\
       &\leq\Lambda^\varepsilon_{3,1s}(\eta) < \Lambda^\varepsilon_{1,2c}(\eta)\leq\Lambda^\varepsilon_{1,2s}(\eta)<\dots,\label{truespectrum}
       \end{align}
Here, $\Lambda^0_{0,k}$, for $k\in\mathbb{N}$,  is a simple eigenvalue,  while $\Lambda^0_{n,k}$, for $n,k\in\mathbb{N}$, is a double eigenvalue of the limit problem in $\Theta$, {\it i.e.} $\Lambda^0_{n,kc}=\Lambda^0_{n,ks}$. We denote by  $\Lambda^0_{n,kc}$ and $\Lambda^0_{n,ks}$ the double eigenvalue corresponding to the cosine and sine eigenfunctions respectively. The explicit expression of $\Lambda^0_{n,k}$ and $U^0_{\Theta, n}$ leads us to compute  the first-order correction term of the asymptotic expansion of $\Lambda^\varepsilon_{n,k}(\eta)$ (see formula \eqref{simpleeig} in the case of a simple eigenvalue $\Lambda^0_{0,k}$ and formula \eqref{multipleeig} in the case of a multiple eigenvalue $\Lambda^0_{n,k}$). This combined with Theorem \ref{thmmain}   provides an asymptotic estimate of the length of spectral bands ${\cal B}^\varepsilon_{n,k}$. We show that the length of ${\cal B}^\varepsilon_{n,k}$ is of order $\varepsilon^{2m}$, for $m\in(0,1/2)$, for any  $k\in\mathbb{N}$ and for $n\in\mathbb{N}\cup\{0\}$ expect for $n=4, 8, 12, \dots$ (see Corollary \ref{cor}). In the latter case, further investigation of higher order terms of the ansatz of $\Lambda^\varepsilon_{n,k}(\eta)$ is to carry out and it is left as an open question to be considered. However, the results loses to be explicit and to make conclusions one needs to provide numeric computations.
 \par
\begin{figure}
\centering
\includegraphics[scale=0.43]{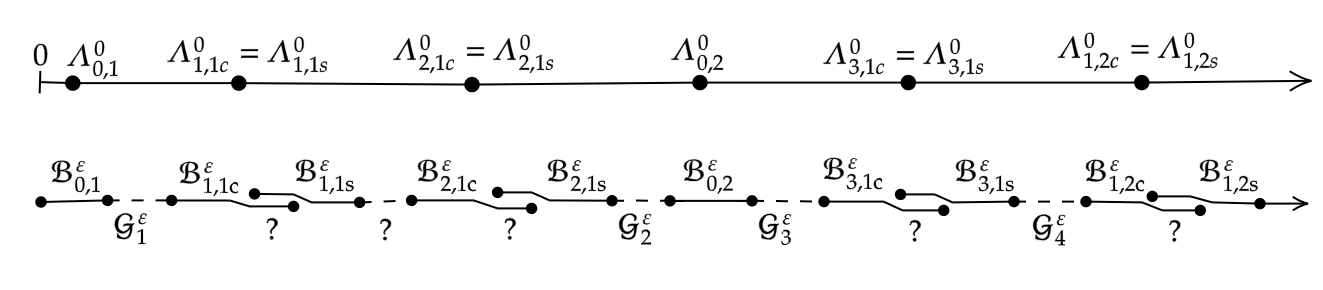}
\caption{\small The spectrum \eqref{spectrumlimprob}  of the limit problem in $\Theta$ and the band-gap structure of the spectrum $\espec^\varepsilon$. Here, for $k\in\mathbb{N}$, $\mathcal{B}^\varepsilon_{0,k}$ denotes the band associated to $\Lambda^\varepsilon_{0,k}(\eta)$, ${\cal B}^\varepsilon_{n,kc}$ and ${\cal B}^\varepsilon_{n,ks}$ are the bands related to $\Lambda^\varepsilon_{n,kc}(\eta)$ and $\Lambda^\varepsilon_{n,ks}(\eta)$, for $n,k\in\mathbb{N}$. The question points $?$ indicate that, up to now, we do not know if the bands ${\cal B}^\varepsilon_{n,kc}$ and ${\cal B}^\varepsilon_{n,ks}$ overlap or are disjoint. Moreover, a further investigation of higher order term is requested to detect the opening of a spectral gap between ${\cal B}^\varepsilon_{1,1s}$ and ${\cal B}^\varepsilon_{2,1c}$.} 
 \label{fig:spectrum}
\end{figure}

In view of the particular position of the zeros of the Bessel function and consequently of the sequence \eqref{truespectrum}, we do not provide a complete result of the existence of spectral gaps ${\cal G}^\varepsilon_n$ (see Figure \ref{fig:spectrum}). More specifically, we only prove the existence of the spectral gaps ${\cal G}^\varepsilon$ between bands generated by eigenvalues $\Lambda^\varepsilon_{n,k}(\eta)$ whose leading term is the simple eigenvalue $\Lambda^0_{0,1}(\eta)$ and the double one $\Lambda^0_{1,1}(\eta)$, the double eigenvalue $\Lambda^0_{2,1}(\eta)$ and the simple one $\Lambda^0_{0,2}(\eta)$,  the simple eigenvalue $\Lambda^0_{0,2}(\eta)$ and the double one $\Lambda^0_{3,1}(\eta)$ and finally the double eigenvalues  $\Lambda^0_{3,1}(\eta)$ and $\Lambda^0_{1,2}(\eta)$ (see Corollary \ref{cor2}).  We highlight that, up to now, we can not detect a gap between the spectral bands generated by eigenvalues whose leading terms are a double  eigenvalues, {\it i.e.} $\Lambda^0_{n,kc}=\Lambda^0_{n,ks}$, for $n,k\in\mathbb{N}$ as well as between the spectral bands ${\cal B}^\varepsilon_{1,1}$ and  ${\cal B}^\varepsilon_{2,1}$ (see Figure \ref{fig:spectrum}) since an analysis of higher order terms of asymptotic expansions of $\Lambda^\varepsilon(\eta)$ is needed and it is left as an open problem.

\bigskip
\subsection{A short review}
The detection of spectral gaps in scalar problems have investigated in \cite{CCNT, HL00, HP01,  N10, Z05}, where the coefficients of the differential operator have high contrast. In our context, the main novelty is that  neighbouring smooth hard fragments of the inhomogeneous plane touch each other at cusp points and the same occurs for the soft fragments which also have geometric irregularities due to cuspidal points. 
In \cite{NRT12} the authors have shown the appearance of the gaps in the spectrum of the Dirichlet and the Neumann problems for the Laplace operator in the plane $\mathbb{R}^2$  perforated by a double-periodic family of circular and isolated holes  located at a small distance and forming thin ligaments. In other words, the interface is smooth. In \cite{FT15} the appearance of the gaps has been done in a more general geometrical settings but for Dirichlet conditions only.
\par In the case of waveguides, the appearance of spectral gaps forbids the wave propagation in the corresponding frequency range. In the literature, there are numerous treatments on the propagation of waves through periodic structures. Two approaches are usually used in order to detect the opening of the gaps: by studying the asymptotic behaviour of eigenvalues of the model problem on a periodicity cell or by seeking for the location of the eigenvalues by means of specific weight estimates, such as the Hardy inequality and the max-min principle. In  \cite{CPNR13, FS08, N10-1} the first approach is used to detect spectral gaps, while in \cite{N10-02, NRT10} the second method is applied. In \cite{NRT12} the authors  have adopted both approaches: the asymptotic method is used for analysing the spectrum of the Dirichlet problem but {\it a priori} estimates of eigenfunctions of the Neumann problem on thin ligaments are necessary to localize the eigenvalues.\\

\smallskip
The paper is organized as follows. In Section $2$ we present the variational formulation of the model problem \eqref{uthetaG}-\eqref{quasipercond2} in the periodicity cell $\perceltheta$ and we discuss the formal asymptotic ans\"{a}tze of the eigenpairs for $m\in (0, 1/2)$. In Section $3$ we discuss the justification procedure. Moreover, we provide the length of the spectral bands ${\cal B}^\varepsilon_{n,k}$ and we detect the opening of the spectral gaps ${\cal G}^\varepsilon$.

\section{Formal asymptotic analysis in the case\\ $\boldsymbol{0<m<1/2}$}
\subsection{The model problem in the periodicity cell}
Let $L^2(\ball)$ and $L^2(\complball)$ be the complex Lebesgue spaces on $\ball$ and $\complball$ endowed with the scalar product $(\cdot, \cdot)_\Theta$ and $(\cdot, \cdot)_{Q\setminus\Theta}$ respectively and let $\perfunc$ be the space of functions in $H^1(\perceltheta)$ satisfying the periodicity conditions \eqref{quasipercond1}. 
The variational setting of the problem \eqref{uthetaG}-\eqref{quasipercond2} reads as
      \begin{align}
      ((\grad+i\eta) \utheta, (\grad+i\eta) V_{\ball})_{\ball} &+\varepsilon^{-1}((\grad+i\eta) \ucomptheta, (\grad+i\eta) V_{\complball})_{\complball}\notag\\
      & = \eiGelf(\eta) \left((\utheta,V_{\ball})_{\ball} + \varepsilon^{-2m}(\ucomptheta, V_{\complball})_{\complball}\right), \label{varformsquare}
       \end{align}
for all $V \in\perfunc$.
 Due to closedness and  positiveness of the sesquilinear form on the left-hand side of \eqref{varformsquare} and thanks to the compactness of the embedding $\perfunc\subset L^2(\perceltheta)$, the  operator $A^\varepsilon(\eta)$ associated to the problem \eqref{varformsquare} is positive, self-adjoint and has a discrete spectrum, given by \eqref{spectrum}. We assume that the eigenfunctions $U^\varepsilon(\cdot,\eta) = (U^{\varepsilon}_{\ball}(\cdot,\eta),U^{\varepsilon}_{\complball,}(\cdot,\eta))$ associated with the identity \eqref{varformsquare} are subject to the orthonormalization conditions
     \begin{equation}
     \label{orthonomcond}
     (U^{\varepsilon}_{\ball,n},U^{\varepsilon}_{\ball,m})_{\ball} + \epsilon^{-2m}(U^{\varepsilon}_{\complball,n}, U^{\varepsilon}_{\complball,m})_{\complball} = \delta_{n,m}, \qquad\mbox{for } n,m\in\mathbb{N},
     \end{equation}
where $\delta_{n,m}$ is the Kronecker symbol. Following similar arguments of \cite{CDN} and in view of  the orthonormalization conditions \eqref{orthonomcond}, we perform the replacement
    \begin{equation*}
    \vtheta(x,\eta) :=\utheta(x,\eta),\qquad \vcomptheta (x,\eta) := \varepsilon^{-m}\ucomptheta(x,\eta).
    \end{equation*}
Hence, the differential equations \eqref{uthetaG}-\eqref{uupsiG} remain invariable as well as the periodicity conditions \eqref{quasipercond1}-\eqref{quasipercond2}, while the transmission conditions become 
    \begin{align}
   \epsilon^{m}\vcomptheta(x,\eta)  &= \vtheta(x, \eta), \hspace{3.4cm} x\in \Gamma, \label{ftcG}\\
    \epsilon^{m-1}\nu\cdot(\grad+i\eta) \vcomptheta(x,\eta) &= \nu\cdot(\grad+i\eta) \vtheta(x,\eta), \hspace{1.5cm} x\in \Gamma.\label{stcG}
    \end{align} 
Thanks to  results obtained in \cite{CDN} where a similar problem for the Laplace operator has been studied for $m\in (0, 1/2)$,  we search for the asymptotic ans\"{a}tze of eigenvalues $\eiGelf$ and eigenfunctions $\{V^\varepsilon_\Theta, V^\varepsilon_{Q\setminus\Theta}\}$ of the form
     \begin{align}
     \eiGelf(\eta) &= \Lambda^0(\eta) + \varepsilon^{2m}\Lambda^1(\eta)+\cdots,\label{ansvalu}\\
     \vtheta(x,\eta) &= V^0_{\ball}(x,\eta) + \varepsilon^{2m}V^1_{\ball}(x,\eta)+\cdots, \hspace{2.1cm}x\in\Theta,\label{ansutheta}\\
     \vcomptheta (x,\eta) &= \varepsilon^{m}V^0_{\complball}(x,\eta) + \varepsilon^{1-m} V^1_{\complball}(x,\eta)+\cdots, \qquad x\in Q\setminus
     \Theta.\label{ansucomptheta}
     \end{align}
Note that the above expansions depend also on the dual variable $\eta\in [-\pi,\pi)^2$.  In order to find the leading and the first-order correction terms, we insert  \eqref{ansvalu}-\eqref{ansucomptheta} into problem \eqref{uthetaG}-\eqref{uupsiG} combined with the new transmission conditions \eqref{ftcG}-\eqref{stcG} and we collect the coefficients of the same powers of $\varepsilon$, obtaining the desired boundary value problems. Note that the feature of the geometry of the inhomogeneous plane enable us to regard the hard and the weakly connected soft   fragments as isolated and independent since they touch at the cusp points only.
\subsection{Problem satisfied by $\boldsymbol{V^0_{\ball}}$}  
The leading term $V^0_{\ball}$ in \eqref{ansutheta} solves the spectral problem
      \begin{align}
      -(\grad+i\eta)^2V^0_{\ball}(x,\eta) &= \Lambda^0 (\eta)V^0_{\ball}(x,\eta), \hspace{0.5cm}x\in\ball,\label{pbvtheta}\\
      V^0_{\ball}(x,\eta) &=0, \hspace{2.5cm}x\in\Gamma.\label{bcvtheta}
      \end{align}
We look for  a solution of the problem \eqref{pbvtheta}-\eqref {bcvtheta} in the form
     \begin{equation*}
     V^0_{\ball}(x,\eta) := e^{-i\eta\cdot x}\mathfrak{V}^0_{\ball}(x, \eta).
     \end{equation*}
Then, $\mathfrak{V}^0_{\ball}(x,\eta)$ is a solution to the problem 
    \begin{align}
    -\laplace\mathfrak{V}^0_{\ball}(x, \eta) &= \Lambda^0 (\eta)\mathfrak{V}^0_{\ball}(x, \eta),\hspace{0.5cm}x\in\ball,\label{specprobl}\\
          \mathfrak{V}^0_{\ball}(x, \eta) &=0.\hspace{2.6cm}x\in\Gamma.\label{Dirchbc}
    \end{align}
The pair $(\Lambda^0
 (\eta), \mathfrak{V}^0_{\ball}(\cdot,\eta))$ is formed by eigenvalue and associated eigenfunction of the Dirichlet Laplacian in the disk $\ball$, hence they are independent of parameter $\eta$.
 In the sequel, we simply write $\mathfrak{V}^0_n$ in place of $\mathfrak{V}^0_{\ball,n}$. \par Thanks to the link between the Bessel functions and eigenpairs of the  Dirichlet Laplacian (see \cite{B58, W95}), we know that the eigenfunctions  $\mathfrak{V}^0_{n,k}$ are given by the Bessel functions $J_{n}$ of the first kind and the eigenvalues $\Lambda^0_{n,k}$ are  the corresponding positive zeros $j_{n,k}$, for $n\in\mathbb{N}\cup\{0\}$ and $k\in\mathbb{N}$. In other words,
         \begin{equation}
         \label{simple}
         \Lambda^0_{0,k}=4j_{0,k}^2, \qquad \mathfrak{V}^0_{0,k}(r)= J_0\left(2j_{0,k}r\right), \quad\mbox{for } k\in\mathbb{N},
         \end{equation}
 and 
          \begin{equation} 
          \label{multiple}
          \Lambda^0_{n,k}=4j_{n,k}^2, \qquad \mathfrak{V}^0_{n,k}(r,\theta)= J_n\left(2j_{n,k}r\right)\left(C_c\cos(n\theta)+C_s\sin(n\theta)\right), \quad\mbox{for } n,k\in\mathbb{N},
          \end{equation}     
 where $C_c$ and $C_s$ are arbitrary constants and $(r, \theta)$ are the polar coordinates. Recall that for fixed $n\in\mathbb{N}\cup\{0\}$, $J_{n}$ has an infinite number of positive real zeros $j_{n,k}$, for $k\in\mathbb{N}$, and any two different Bessel functions $J_{n}$ and $J_{l}$ do not get common roots except for $j_{n,0}=j_{l,0}=0$, for $n,l\in\mathbb{N}$ (see \cite{BK78}). Therefore, the spectrum of the problem \eqref{pbvtheta}-\eqref{bcvtheta}, being independent of the Floquet parameter $\eta$,  consists of the sequence \eqref{spectrumlimprob}.
 \vspace{0.5cm}

 \subsection{Problem satisfied by $\boldsymbol{V^0_{\complball}}$} 
 \begin{figure}
 \centering
 \begin{minipage}{0.3\textwidth}
 \centering
 \includegraphics[scale=0.22]{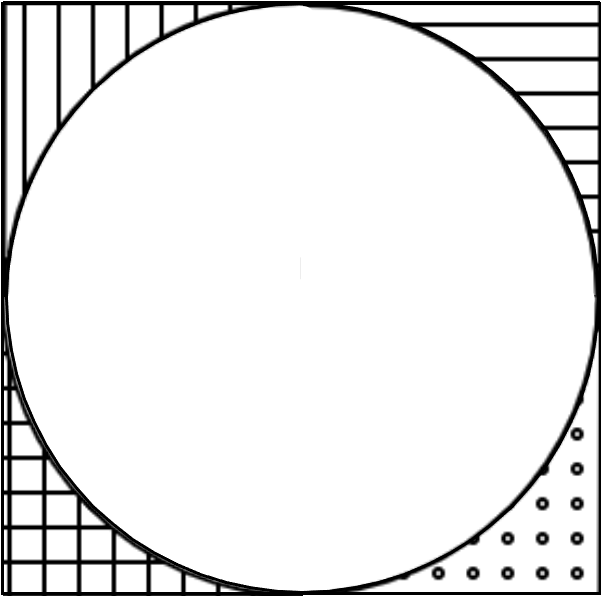}
 \caption*{$\perceltheta$}
 \end{minipage} %
 \begin{minipage}{0.1\textwidth}
 \centering
 \includegraphics[scale=0.4]{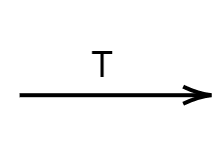}
 \end{minipage} %
 \hspace{0.5cm}
 \begin{minipage}{0.3\textwidth}
 \centering
 \includegraphics[scale=0.22]{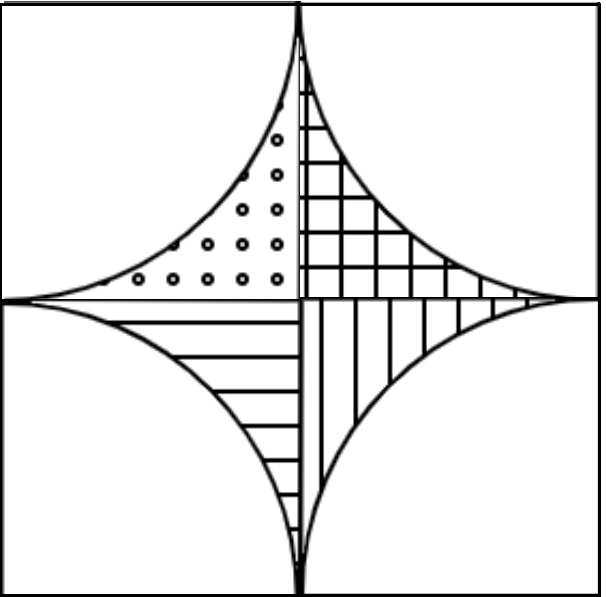}
 \caption*{$\percelupsilon$}
 \end{minipage}
 \caption{\small It shows how the map $T$ trasforms the periodicity cell $\perceltheta$ into that $\percelupsilon$.}
 \label{figure:mapT}
 \end{figure}
 
 The leading term $V^0_{\complball}$ of the expansion \eqref{ansucomptheta} is a solution of the problem
       \begin{align*}
        -(\grad + i\eta)^2V^0_{\complball}(x,\eta) &=0, \hspace{0.9cm}x\in\complball,\\
        \nu\cdot(\grad+i\eta)V^0_{\complball}(x,\eta) &=0, \hspace{0.8cm}x\in\Gamma.
       \end{align*}
 We look for a solution of the form 
      \begin{equation*}
      V^0_{Q\setminus\Theta}(x,\eta): = e^{-i\eta\cdot x}\mathfrak{V}^0_{Q\setminus\Theta}(x, \eta).
      \end{equation*}
 The function $\mathfrak{V}^0_{Q\setminus\Theta}$ satisfies the boundary value problem 
      \begin{align}
      -\laplace\mathfrak{V}^0_{Q\setminus\Theta}(x, \eta) &=0,\qquad x\in Q\setminus\Theta,\label{pbB0}\\
      \partial_\nu \mathfrak{V}^0_{Q\setminus\Theta}(x,\eta)&=0,\qquad x\in \Gamma.\label{bcB0}
      \end{align}
  Note that the set $Q\setminus\Theta$ is disconnected in $\perceltheta$. However, thanks to the periodicity conditions \eqref{quasipercond1}, we can switch periodicity cell from $\perceltheta$ to $\percelupsilon$ such that the set $Q\setminus\Theta$ turns into a connected set $\Upsilon$ in $\percelupsilon$ (see Figure \ref{figure:mapT}). Hence, solving a boundary value problem  in $Q\setminus\Theta$ together with periodicity conditions is equivalent to solve the same boundary value problem in the connected set $\Upsilon$ without the periodicity conditions. This allows us to deal with problems  in the periodicity cell $\percelupsilon$ using the classical Sobolev spaces on the connected set $\Upsilon$. \par Let $Q^j_\Theta$ and $Q^j_\Upsilon$ be defined by
        \begin{align*}
        Q^1_k &:=\{x=(x_1,x_2)\in\omega_k\hspace{0.02cm}:\hspace{0.02cm} 0<x_j\leq 1/2, \hspace{0.05cm}j=1,2 \},\\
        Q^2_k &:=\{x=(x_1,x_2)\in\omega_k\hspace{0.02cm}:\hspace{0.02cm} -1/2\leq x_1<0,\hspace{0.1cm}0<x_2\leq 1/2 \},\\
        Q^3_k &:=\{x=(x_1,x_2)\in\omega_k\hspace{0.02cm}:\hspace{0.02cm} -1/2\leq x_j<0, \hspace{0.05cm}j=1,2 \},\\
        Q^4_k &:=\{x=(x_1,x_2)\in\omega_k\hspace{0.02cm}:\hspace{0.02cm} 0<x_1\leq 1/2,\hspace{0.1cm}-1/2\leq x_2<0 \},
        \end{align*}
 for $k=\Theta, \Upsilon$. We define the map   $\hat{x}=T(x)$ where ${\rm T}\hspace{0.01cm}:\hspace{0.01cm}\perceltheta\rightarrow\percelupsilon  $   is defined by
      \begin{align}
      (x_1,x_2)\in Q^1_{\Theta}&\mapsto {\rm T}(x_1,x_2):= (x_1-1/2, x_2-1/2) \in Q^3_{\Upsilon},\label{defT1}\\
      (x_1,x_2)\in Q^2_{\Theta} &\mapsto {\rm T}(x_1,x_2):= (x_1+1/2, x_2-1/2)\in Q^4_{\Upsilon},\notag\\
      (x_1,x_2)\in Q^3_{\Theta}  &\mapsto {\rm T}(x_1,x_2):= (x_1+1/2, x_2+1/2)\in Q^1_{\Upsilon},\notag\\
      (x_1,x_2)\in Q^4_{\Theta} &\mapsto {\rm T}(x_1,x_2):= (x_1-1/2, x_2+1/2)\in Q^2_{\Upsilon}.\label{defT2}
      \end{align}
This implies that ${\rm T}$ maps the function $V^0_{Q\setminus\Theta}(x, \eta)$, for $x\in Q\setminus\Theta$, in the function $V^0_{\Upsilon}(\hat{x}, \eta)$, for $\hat{x}\in\Upsilon$. Now, we look for 
    \begin{equation}
    \label{Vupsilon}
    V^0_{\Upsilon}(\hat{x}, \eta) := e^{-i\eta\cdot\hat{x}} \mathfrak{B}^0_{\Upsilon}(\hat{x}, \eta).
    \end{equation}
The function $\mathfrak{B}^0_{\Upsilon}(\hat{x}, \eta)$ satisfies the same problem as $\mathfrak{B}^0_{Q\setminus\Theta}(x, \eta)$ but on different region, \textit{i.e.}
     \begin{align*}
     -\laplace \mathfrak{B}^0_{\Upsilon}(\hat{x}, \eta) &=0, \qquad \hat{x}\in\Upsilon,\\
     \partial_v \mathfrak{B}^0_{\Upsilon}(\hat{x}, \eta)&=0, \qquad \hat{x}\in \Gamma_\Upsilon,
     \end{align*}
with $\Gamma_{\Upsilon} := \partial\Upsilon$.
Hence, due to the connectedness of $\Upsilon$, $\mathfrak{B}^0_{\Upsilon}(\hat{x}, \eta)$ is a constant function with respect to the variables $\hat{x}$, \textit{i.e.}
    \begin{equation}
    \label{constantc0}
    \mathfrak{B}^0_{\Upsilon}(\hat{x}, \eta) = c^0(\eta).
    \end{equation}
In light of the definition \eqref{defT1}-\eqref{defT2} of the map $T$  and due  to \eqref{Vupsilon}, we get that
   \begin{align}
   \label{Vinverse}
   V_{Q\setminus\Theta}^0(x, \eta)=V_{\Upsilon}^0(\hat{x}, \eta) = e^{-i\eta\cdot\hat{x}} \mathfrak{B}^0_{\Upsilon}(\hat{x}, \eta) =e^{-i\eta\cdot\hat{x}} c^0(\eta). 
   \end{align}
In view of \eqref{defT1}-\eqref{defT2}, it is easy to check that 
    \begin{equation*}
    e^{-i\eta\cdot\hat{x}} = e^{-i\eta\cdot x}g_x(\eta),
    \end{equation*}
where $g_x(\eta)$ is defined by 
      \begin{equation}
      \label{mapg}
      g_x(\eta) :=
        \begin{cases}
        e^{i(\eta_1/2+\eta_2/2)},& x\in Q^1_{\Theta},\\
        e^{i(-\eta_1/2+\eta_2/2)},& x\in Q^2_{\Theta},\\
        e^{-i(\eta_1/2+\eta_2/2)},& x\in Q^3_{\Theta},\\
        e^{i(\eta_1/2-\eta_2/2)},& x\in Q^4_{\Theta}.\\
        \end{cases}
      \end{equation}
This combined with \eqref{Vinverse} implies that 
     \begin{equation*}
     V^0_{Q\setminus\Theta}(x,\eta) = g_x(\eta) e^{-\eta\cdot x} c^0(\eta).
     \end{equation*} 
     
\subsection{Problem satisfied by $\boldsymbol{V^1_{\complball}}$} The first-order correction term $V^1_{\complball}$ in \eqref{ansucomptheta} is a the solution to the problem
     \begin{align*}
     -(\grad+i\eta)^2V^1_{\complball}(x,\eta) &= \Lambda^0 V^0_{\complball}(x, \eta),\hspace{1.5cm} x\in\complball,\\
     \nu\cdot(\grad+i\eta)V^1_{\complball} (x, \eta) &= \nu\cdot(\grad+i\eta)V^0_{\ball} (x), \qquad x\in\Gamma.
     \end{align*}
Thanks to the map $T$ given by \eqref{defT1}-\eqref{defT2}, we know that $V^1_{Q\setminus\Theta}(x, \eta) = V^1_{\Upsilon}(\hat{x}, \eta)$. Hence, we look for a solution $V^1_\Upsilon$ of the form 
    \begin{equation*}
    V^1_\Upsilon(\hat{x}, \eta) = e^{-i\eta\cdot\hat{x}} \mathfrak{V}^1_\Upsilon(\hat{x}, \eta).
    \end{equation*}
The function  $\mathfrak{V}^1_{\Upsilon}(\hat{x}, \eta)$ satisfies the boundary value problem 
      \begin{align}
      -\laplace \mathfrak{V}^1_{\Upsilon} (\hat{x},\eta)&= \Lambda^0_{n,k}c^0(\eta), \hspace{1.7cm}\hat{x}\in\Upsilon,\label{B1ups}\\
      \normalder \mathfrak{V}^1_{\Upsilon}(\hat{x},\eta) &= g^{-1}_x(\eta)\normalder  \mathfrak{V}^0_{n,k}(\hat{x}), \hspace{0.5cm}\hat{x}\in\Gamma_{\Upsilon},\label{B1upsbc}
      \end{align}
for $n\in\mathbb{N}\cup\{0\}$ and for $k\in\mathbb{N}$, where $g^{-1}$ is the inverse of the map $g$ defined by \eqref{mapg} and $c^0(\eta)$ is given by \eqref{constantc0}. 
The compatibility condition to the problem \eqref{B1ups}-\eqref{B1upsbc} reads as
     \begin{equation*}
     \int_{\Upsilon} \Lambda^0_{n,k} c^0(\eta) dx = -  \int_{\Gamma_\Upsilon}g^{-1}_x(\eta)\normalder\mathfrak{V}^0_{n,k}(\hat{x})ds_{\hat{x}}, \qquad\mbox{for }n\in\mathbb{N}\cup\{0\}, k\in\mathbb{N},
     \end{equation*}
which implies that the constant $c^0(\eta)$ is given by 
    \begin{equation}
    \label{constant}
    c^0(\eta) 
    := -\frac{1}{\Lambda^0_{n,k}(1-\pi/4)} \int_{\Gamma_\Upsilon}g^{-1}_x(\eta)\normalder\mathfrak{V}^0_{n,k}(\hat{x})ds_{\hat{x}}, \quad\mbox{for }n\in\mathbb{N}\cup\{0\}, k\in\mathbb{N}.
    \end{equation}
    
\subsubsection{Simple eigenvalues}
We assume that $n=0$. Recall that the derivative of the Bessel function $J_{0}$ is given by the formula (see \textit{e.g.} \cite[Chapter VI]{B58}) 
     \begin{equation*}
     \frac{d}{dx}J_0(x) = -J_1(x).
     \end{equation*}
In view of the definition of the function \eqref{mapg} and due to the formulas \eqref{simple}, from the equality \eqref{constant} we deduce that
      \begin{align}
      c^0_{0,k}(\eta) &= -\frac{1}{\Lambda^0_{n,k}(1-\pi/4)}\frac{d}{dr}J_0\left(2j_{0,k}r \right)_{\big|r=1/2}\biggl(\int_{0}^{\pi/2} e^{i(\eta_1/2+\eta_2/2)}  d\theta  + \int_{\pi/2}^{\pi} e^{-i(\eta_1/2-\eta_2/2)} d\theta   \notag\\
      &\hspace{1cm}+  \int_{\pi}^{3\pi/2} e^{-i(\eta_1/2+\eta_2/2)}  d\theta +\int_{3\pi/2}^{2\pi} e^{i(\eta_1/2-\eta_2/2)}   d\theta\biggr)\notag\\
      &= \frac{1}{\Lambda^0_{n,k}(1-\pi/4)}\pi j_{0,k} J_1\left(j_{0,k}\right)( e^{i(\eta_1/2+\eta_2/2)} +   e^{-i(\eta_1/2-\eta_2/2)} + e^{-i(\eta_1/2+\eta_2/2)}+e^{i(\eta_1/2-\eta_2/2)}   )\notag\\
      &= \frac{\pi}{\Lambda^0_{n,k}(1-\pi/4)}j_{0,k}\left[e^{i\eta_1/2}+e^{-i\eta_1/2}\right]\left[e^{i\eta_2/2}+e^{-i\eta_2/2}\right] J_1\left(j_{0,k} \right)\notag\\
      &=\frac{\pi}{j_{0,k}(1-\pi/4)}J_1\left(j_{0,k}\right)\cos\left(\frac{\eta_1}{2}   \right) \cos\left(\frac{\eta_2}{2}\right) \label{csimple}, \qquad\mbox{for } k\in\mathbb{N}.
      \end{align}

\subsubsection{Multiple eigenvalues}
Now, assume that $n\in\mathbb{N}$. Recall that the derivative of the Bessel function $J_n(x)$ is given by the recurrence formula (see \textit{e.g.} \cite[Chapter VI]{B58}) 
    \begin{align*}
    \frac{d}{dx}J_n(x) &= \frac{1}{2}\left( J_{n-1}(x)-J_{n+1}(x)  \right).
    \end{align*}
Then, in view of the definition \eqref{multiple}, a direct computation yields  
    \begin{align}
    c^0_{n,k}(\eta) &= -\frac{1}{\Lambda^0_{n,k}(1-\pi/4)}\frac{d}{dr}J_n\left(2j_{n,k}r \right)_{\big|r=1/2}\biggl(\int_{0}^{\pi/2} e^{i(\eta_1/2+\eta_2/2)}(C_c\cos(n\theta) 
     + C_s\sin(n\theta))  d\theta \notag\\
     &\quad+ \int_{\pi/2}^{\pi} e^{-i(\eta_1/2-\eta_2/2)} (C_c\cos(n\theta) + C_s\sin(n\theta)) d\theta \notag \\
     &\quad +  \int_{\pi}^{3\pi/2} e^{-i(\eta_1/2+\eta_2/2)} (C_c\cos(n\theta) + C_s\sin(n\theta)) d\theta \notag\\
    &\quad +\int_{3\pi/2}^{2\pi} e^{i(\eta_1/2-\eta_2/2)} (C_c\cos(n\theta) + C_s\sin(n\theta)) d\theta\biggr) \notag\\
    &=
       \begin{cases}
       0, \hspace{1.9cm} n=4,8,12,\dots,\\
      \alpha_1, \hspace{1.7cm} n=2,6,10,\dots,\\
      \alpha_2 , \hspace{1.7cm} n=1,5,9,\dots,\\
       \alpha_3  , \hspace{1.7cm} n=3,7,11,\dots, \label{cmultiple}
       \end{cases}
    \end{align}
where
    \begin{align*}
    \alpha_1 &:=-\frac{4C_s}{nj_{n,k}(1-\pi/4) } \left(J_{n-1}\left(j_{n,k}\right) - J_{n+1}\left(j_{n,k}\right) \right) \sin ( \frac{\eta_1}{2} )\sin ( \frac{\eta_2}{2} ),\\
    \alpha_2 &:=  -\frac{2i}{nj_{n,k}(1-\pi/4) } \left(J_{n-1}\left(j_{n,k}\right) - J_{n+1}\left(j_{n,k}\right) \right) \left(C_c\sin ( \frac{\eta_1}{2} )\cos ( \frac{\eta_2}{2} )+C_s\cos ( \frac{\eta_1}{2} )\sin ( \frac{\eta_2}{2} )\right),\\
    \alpha_3 &:=  \frac{2i}{nj_{n,k}(1-\pi/4) } \left(J_{n-1}\left(j_{n,k}\right) - J_{n+1}\left(j_{n,k}\right) \right) \left(C_c\sin ( \frac{\eta_1}{2} )\cos ( \frac{\eta_2}{2} )-C_s\cos ( \frac{\eta_1}{2} )\sin ( \frac{\eta_2}{2} )\right).
    \end{align*}
\subsection{Problem satisfied by $\boldsymbol{V^1_{\ball}}$} 
The first-order correction term $V^1_{\ball}$ verifies the problem
     \begin{align*}
     -(\grad+i\eta)^2V^1_{\ball}(x,\eta) -\Lambda^0V^1_{\ball}(x,\eta) &= \Lambda^1(\eta)V^0_{\ball}(x), \hspace{1.1cm} x\in\ball,\\
     V^1_{\ball}(x,\eta)&=V^0_{\complball}(x,\eta), \hspace{1.3cm}x\in\Gamma.
     \end{align*}
We look for $V^1_{\ball}$   in the form $$V^1_{\ball} (x,\eta):= e^{-i\eta\cdot x}\mathfrak{V}^1_{\ball}(x,\eta),$$ where $\mathfrak{V}^1_{\ball}$ is a solution to
    \begin{align}
    -\laplace \mathfrak{V}^1_{\ball} (x,\eta) - \Lambda^0_n\mathfrak{V}^1_{\ball}(x,\eta) &= \Lambda^1_n(\eta)\mathfrak{V}^0_{n,k}(x), \qquad x\in\ball,\label{pbV'}\\
    \mathfrak{V}^1_{\ball}(x,\eta) &=g_x(\eta)c^{0}_{n,k}(\eta), \hspace{1cm} x\in \Gamma,\label{bcV'}
    \end{align}  
where $c^{0}_{n,k}(\eta)$ is given by formula \eqref{csimple} if $n=0$ and \eqref{cmultiple} if $n\in\mathbb{N}$.\par Recall that $\mathfrak{V}^1_{\ball}(x,\eta) = g_x(\eta)c_{0,k}(\eta)$ on $\Gamma$.  In the case of simple eigenvalues $\Lambda^0_{0,k}$, for $k\in\mathbb{N}$, the Fredholm alternative leads us to a single compatibility condition 
           \begin{align}
          \Lambda^1_{0,k}(\eta) (\mathfrak{V}^0_{0,k}, \mathfrak{V}^0_{0,k})_{\Theta} &= \int_{\Gamma}
            g_x(\eta)c_{0,k}(\eta)  \frac{d}{dr} \mathfrak{V}^0_{0,k}(r) d\theta\notag\\
            &=\frac{d}{dr} J_{0,k}(2j_{0,k}r)_{\big|r=1/2}\int_{\Gamma}g_x(\eta)c^0_{0,k}(\eta)d\theta\notag\\
            &=2j_{0,k} J_{1,k}(j_{0,k})c^0_{0,k}(\eta)\frac{\pi}{2}4\cos\left( \frac{\eta_1}{2}\right)\cos\left( \frac{\eta_2}{2}\right)\notag\\
            &=  \frac{2\pi}{1-\pi/4}\left(J_{1,k}(j_{0,k})  \cos\left( \frac{\eta_1}{2}\right)\cos\left( \frac{\eta_2}{2}\right) \right)^2.\label{simpleeig}
            \end{align}
Assume, now, that $\Lambda^0_{n,k}$ is an eigenvalue with multiplicity two. For convenience, we denote the corresponding eigenfunctions  by
      \begin{equation*}
      \mathfrak{V}^0_{n, kc}(r, \theta) := J_n\left(2j_{n,k} r \right)\cos(n\theta), \qquad \mathfrak{V}^0_{n, ks}(r, \theta) := J_n\left(2j_{n,k} r \right)\sin(n\theta).
      \end{equation*}
Hence, we predict that the term $\mathfrak{V}^0_{n,k}(r,\theta)$ takes the form
       \begin{equation*}
       \mathcal{V}^0_{n,k\hbar} (r,\theta): = a^\hbar_c\mathfrak{V}^0_{n,kc}(r, \theta) + a^\hbar_{s}\mathfrak{V}^0_{n,ks}(r, \theta), \qquad\mbox{for } \hbar=c,s,
       \end{equation*}
\textit{i.e.} it is a linear combination of the eigenfunctions $\mathfrak{V}^0_{n,kc}$ and $\mathfrak{V}^0_{n,ks}$.
We require that the vector $a^\hbar=(a^\hbar_c,a^\hbar_{s})\in\mathbb{C}^2$ satisfies the orthonormalization condition
     \begin{equation*}
     (a^\hbar,a^\ell) =a^\hbar_c\overline{a^\ell_c}+a^\hbar_s\overline{a^\ell_s} = \delta_{\hbar \ell}, \qquad\mbox{for } \hbar, \ell=c,s.
     \end{equation*}
Therefore,  $\mathcal{V}^1_{\ball} = \mathcal{V}^1_{n,k\hbar} $ verifies the problem 
    \begin{align*}
    -\laplace\mathcal{V}^1_{n,k\hbar}(x,\eta)-\Lambda^0_n\mathcal{V}^1_{n,k\hbar}(x,\eta) &= \Lambda^1_{n, k\hbar}\mathcal{V}^0_{n,k\hbar}(x,\eta), \qquad x\in\ball,\\
    \mathcal{V}^1_{n,k\hbar}(x,\eta) &= g_x(\eta)c^0_{n,k}(\eta), \hspace{1.3cm} x\in\Gamma.
    \end{align*}
Then, the Fredholm alternative leads us to the two compatibility conditions
    \begin{align*}
    \Lambda^1_{n, k\hbar}(\eta) (\mathcal{V}^1_{n,k\hbar}(x,\eta), \mathfrak{V}^0_{n,k\ell} )&= (\partial \mathfrak{V}^0_{n, k\ell}, \mathcal{V}^1_{n,k\hbar})_{\Gamma}, \qquad\mbox{for }  \ell=c,s.
    \end{align*}
In the algebraic form, $\Lambda^1_{n,kc}$ and $\Lambda^1_{n,ks}$ are eigenvalues with corresponding eigenvectors $a^c$ and $a^s$ of the matrix
         \begin{align*}
         M &:= \frac{1}{j_{n,k}(1-\pi/4)}\left(J_{n-1}\left(j_{n,k}\right) - J_{n+1}\left(j_{n,k}\right) \right) \\
         &\quad\times
         \begin{pmatrix}
         \left(\int_{\Gamma} g_x(\eta) \cos (n\theta )d\theta  \right)^2 &\int_{\Gamma} g_x(\eta) \cos (n\theta )d\theta  \int_{\Gamma} g_x(\eta) \sin (n\theta )d\theta \\
         \int_{\Gamma} g_x(\eta) \cos (n\theta )d\theta  \int_{\Gamma} g_x(\eta) \sin (n\theta )d\theta  & \left(\int_{\Gamma} g_x(\eta) \sin (n\theta )d\theta  \right)^2
         \end{pmatrix}.
         \end{align*}
Therefore, it is easy to check that the eigenvalues of $M$ are given by 
       \begin{equation*}
         \Lambda^1_{n,kc}(\eta) =0 ,\qquad \Lambda^1_{n,ks}(\eta) ={\rm tr}(M),
       \end{equation*}
where the trace of $M$ is 
        \begin{align}        
        \Lambda^1_{n, ks}(\eta)&:= \frac{1}{j_{n,k}(1-\pi/4)}\left(J_{n-1}\left(j_{n,k}\right) - J_{n+1}\left(j_{n,k}\right) \right)\notag\\
        &\quad\times 
            \begin{cases}
            0, & n=4,8,12,\dots\\
            \frac{64}{n^2}\sin^2\left( \frac{\eta_1}{2}\right)\sin^2\left( \frac{\eta_2}{2}\right),&  n= 2,6,10,\dots\\
            -\frac{16}{n^2}\left(\sin^2\left( \frac{\eta_1}{2}\right)\cos^2\left( \frac{\eta_2}{2}\right) +\cos^2\left( \frac{\eta_1}{2}\right)\sin^2\left( \frac{\eta_2}{2}\right)   \right),&  n\quad\text{odd}.\label{multipleeig}
            \end{cases}
        \end{align}

\section{Asymptotic structure of the spectrum}
\subsection{Justification}
\begin{figure}
\centering
\begin{minipage}{0.3\textwidth}
\centering
\includegraphics[scale=0.50]{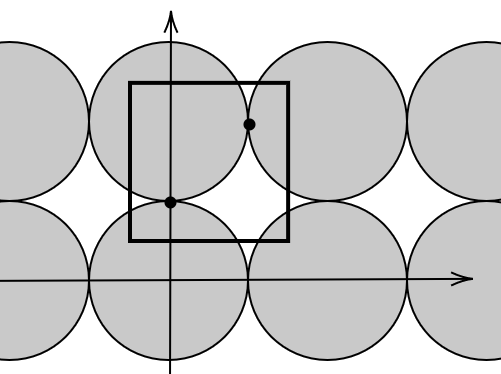}
\end{minipage} 
\hspace{2.4cm}
\begin{minipage}{0.3\textwidth}
\centering
\includegraphics[scale=0.80]{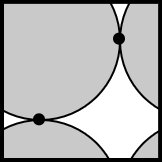}
\caption*{$\omega$}
\end{minipage}
\caption{\small The new choice of the periodicity cell $\omega$}
\label{figure:cellc}
\end{figure}
In this section we comment on the justification of the constructed  formal asymptotics. Our aim is to exploit the results obtained in \cite{CDN}, where a similar problem for the Laplace operator is dealt with. To reach this goal,
we make some changes to our previous analysis: a new version of periodicity cell must be introduced and the Laplace operator and the \textquotedblleft pure\textquotedblright\, normal derivative are required in the model problem to handle the same problem involved in \cite{CDN}. In order to recover the geometry adopted in \cite{CDN}, we choose an alternative periodicity cell such that all the cusps lie in the interior of the cell. More specifically, we translate the unit square $Q=(-1/2, 1/2)^2$ of vector $(1/4, 3/4)$ obtaining $Q'= (-1/4, 3/4)\times (1/4, 5/4)$. Then, the new periodicity cell $\omega$ is defined by $\Xi \cup (\omega\setminus\Xi)$, where $\Xi:= Q'\cap\Omega_1$ (see Figure \ref{figure:cellc}).  
Choosing $\omega$ as a periodicity cell and applying the equivalent version of the Gelfand transform \eqref{newGelftransint} to the problem \eqref{utheta}-\eqref{sbc}, we obtain that  the model problem \eqref{uthetaG}-\eqref{sbcG} along with the periodicity conditions \eqref{quasipercond1}-\eqref{quasipercond2}turns into the model problem \eqref{modlpb}-\eqref{modpb1} together with the quasi-periodicity conditions \eqref{qptrue1}-\eqref{qptrue3}. The integral identity of \eqref{modlpb}-\eqref{qptrue3} reads as   
      \begin{align}
      (\grad  G^\varepsilon_{\Xi}, \grad\psi_\Xi)_{\Xi}&-\varepsilon^{-1}(\grad  G^\varepsilon_{\omega\setminus\Xi}, \grad\psi_{\omega\setminus\Xi})_{\omega\setminus\Xi} \notag\\ 
      &=\Lambda^\varepsilon(\eta)\left((G^\varepsilon_{\Xi}, \psi_\Xi)_{\Xi} + \varepsilon^{-2m}( G^\varepsilon_{\omega\setminus\Xi}, \psi_{\omega\setminus\Xi})_{\omega\setminus\Xi}\right) ,\label{secvarform}
      \end{align}
for any $\psi\in H^1_{\rm qp}$, where $H^1_{\rm qp}$ denotes the subspaces of the Sobolev space $H^1(\omega)$ satisfying the quasi-periodicity conditions \eqref{qptrue1} for $\eta\in [-\pi, \pi)^2$. Since the sesquilinear form on the left-hand side of \eqref{secvarform} is closed and positive and due to compactness of the embedding $H^1_{\rm qp}\subset L^2(\omega)$, the operator $A^\varepsilon_{\rm qp}(\eta)$ associated to \eqref{secvarform} is positive, self-adjoint and its spectrum  is discrete and consists of the sequence \eqref{spectrum}. We assume also that the eigenfunctions $\{G^\varepsilon_\Xi(\cdot, \eta), G^\varepsilon_{\omega\setminus\Xi}(\cdot, \eta) \}$ are subject to the orthonormalization conditions \eqref{orthonomcond}.
\par  
We predict that the formal asymptotic expansions of the eigenpairs $(\Lambda^\varepsilon, \{G^\varepsilon_\Xi, G^\varepsilon_{\omega\setminus\Xi}\})$ take the form \eqref{ansvalu}-\eqref{ansucomptheta}. However, the boundary value problems satisfied by the terms involved in the ans\"{a}tze are different. Indeed, since the two versions of the Gelfand transform \eqref{Gelfandtran} and \eqref{newGelftransint} are linked by the relationship 
   \begin{equation*}
   U(x,\eta) = e^{-i\eta\cdot x}G(x,\eta),
   \end{equation*}
we deduce that the leading and the first-correction terms of the expansions of $G^\varepsilon_\Xi$ satisfy the problems \eqref{specprobl}-\eqref{Dirchbc} and \eqref{pbV'}-\eqref{bcV'}, while the leading term and the first-corrector of the ans\"{a}tze of $G^\varepsilon_{\omega\setminus\Xi}$ are solutions of the boundary value problems  \eqref{pbB0}-\eqref{bcB0} and \eqref{B1ups}-\eqref{B1upsbc}. 
\par The Floquet parameter $\eta\in[-\pi, \pi)^2$ does not represent a trouble due to the continuous dependence of the spectrum \eqref{truespectrum} on  $\eta$ (see \cite{K95}).
This combined with the transform \eqref{newGelftransint} and the change of the periodicity cell,  enables us to repeat the same arguments as  \cite[Theorem 3.1]{CDN} to justify \eqref{ansvalu}-\eqref{ansucomptheta}. 
For reader's convenience, we just state the claim of \cite[Theorem 3.1]{CDN}  in our setting.
  \begin{thm}
   \label{thmmain}
   For $m\in (0,1/2) $ and for any $n,k\in\mathbb{N}$, there exist  $\varepsilon_{n,k}>0$ and  $C_{n,k}>0$ such that for any parameter $\eta\in [-\pi, \pi)^2$ of the Gelfand transform \eqref{Gelfandtran}, the eigenvalues $\Lambda^\varepsilon_{n,k}$ of the problem \eqref{uthetaG}-\eqref{sbcG} along with the periodicity conditions \eqref{quasipercond1}-\eqref{quasipercond2} in the periodicity cell $\perceltheta$ and the eigenvalues  $\Lambda^0_{n,k}$ of the limit problem \eqref{pbvtheta}-\eqref{bcvtheta} are related as follows
           \begin{equation}
            \label{justasy}
             |\Lambda^\varepsilon_{n,k}(\eta)-\Lambda^0_{n,k} - \varepsilon^{2m} \Lambda^1_{n,k}(\eta) | \leq   C_{n,k} \varepsilon^{\gamma},\qquad \varepsilon\in (0,\varepsilon_{n,k}),
           \end{equation}
   with $\gamma=\min\{3m,1 \}$ and $\displaystyle C_{n,k} := \max_{\eta\in [-\pi, \pi)^2} C_{n,k}(\eta)$.
   \end{thm}

\subsection{On the bands and the gaps}
In this subsection, we provide an asymptotic estimate of the length of the spectral bands ${\cal B}^\varepsilon_{n,k}$ and we discuss the opening of the spectral gaps ${\cal G}^\varepsilon$.\par
From  Theorem \ref{thmmain} and formulas \eqref{simpleeig} and \eqref{multipleeig}, it follows the next corollary about the length of spectral bands.
    \begin{cor}
     \label{cor}
     For $k\in\mathbb{N}$, the lengths $L^\varepsilon_{n,k}$ of the spectral bands are given by 
             \begin{align*}
             L^\varepsilon_{0,k} =\varepsilon^{2m} \frac{2\pi}{1-\pi/4}J^2_{1}(j_{0,k}) + O(\varepsilon^\gamma), &\\
             L^\varepsilon_{n,k} =\varepsilon^{2m} \frac{64}{j_{n,k}n^2(1-\pi/4)}(J_{n-1}(j_{n,k})-J_{n+1}(j_{n,k}) )+O(\varepsilon^\gamma),  &\hspace{0.5cm}\mbox{for}\hspace{0.1cm} n=2,6, 10,\dots,\\
             L^\varepsilon_{n,k} = \varepsilon^{2m} \frac{16}{j_{n,k}n^2(1-\pi/4)}(J_{n-1}(j_{n,k})-J_{n+1}(j_{n,k}) )+O(\varepsilon^\gamma), &\hspace{0.5cm} \mbox{for}\hspace{0.1cm} n=\mbox{odd},\\
             L^\varepsilon_{n,k} = O(\varepsilon^{2m}), &\hspace{0.5cm}\mbox{for}\hspace{0.1cm} n =4,8,12,\dots,
             \end{align*}
     where $\gamma=\min\{3m,1 \}$.
     \end{cor}
\begin{proof}
From \eqref{simpleeig} and \eqref{multipleeig}, we have 
 \begin{align*}
    & \Lambda^0_{0,k}-C_{0,k} \varepsilon^{\gamma}\leq \Lambda^\varepsilon_{0,k}(\eta)\leq \Lambda^0_{0,k} + \varepsilon^{2m}\frac{2\pi}{1-\pi/4}J^2_{1}(j_{0,k})+C_{0,k} \varepsilon^{\gamma},\\ 
     &\Lambda^0_{n,k}-C_{n,k} \varepsilon^{\gamma}\leq \Lambda^\varepsilon_{n,k}(\eta)\leq \Lambda^0_{n,k} +  \varepsilon^{2m}\left( \frac{64}{j_{n,k}n^2(1-\pi/4)}(J_{n-1}(j_{n,k})-J_{n+1}(j_{n,k}) ) \right)+C_{n,k} \varepsilon^{\gamma},\\
     &\hspace{12.5cm}\mbox{for}\hspace{0.1cm} n=2,6, 10,\dots,\\
     &\Lambda^0_{n,k}- \varepsilon^{2m}\left( \frac{16}{j_{n,k}n^2(1-\pi/4)}(J_{n-1}(j_{n,k})-J_{n+1}(j_{n,k}) )\right)-C_{n,k} \varepsilon^{\gamma}\leq \Lambda^\varepsilon_{n,k}(\eta)\leq \Lambda^0_{n,k}+C_{n,k} \varepsilon^{\gamma},\\ &\hspace{12.5cm}\mbox{for}\hspace{0.1cm} n=\mbox{odd},
    \end{align*}
which implies the lengths of spectral bands.
\end{proof}
Note that the lengths of the spectral bands ${\cal B}^\varepsilon_{n,k}$ for $n=4, 8, \dots$ are not determined  because of \eqref{multipleeig} and further computations of higher order terms in  the ans\"{a}tze of the eigenpairs $(\Lambda^\varepsilon, \{U^\varepsilon_\Theta, U^\varepsilon_{Q\setminus\Theta}\})$ are necessary.\par
Now, let us investigate the opening of the spectral gaps ${\cal G}^\varepsilon$ in the band-gap structure of the spectrum  \eqref{truespectrum} of the problem \eqref{uthetaG}-\eqref{quasipercond2}. Since the spectrum \eqref{spectrumlimprob} is related to the zeros of the Bessel functions $J_n$, we can not give a complete result about the existence of the spectral gaps. More specifically, we detect the spectral gaps ${\cal G}^\varepsilon$ only in the following cases:  between bands generated by eigenvalues $\Lambda^\varepsilon_{n,k}(\eta)$ whose leading term is given by the simple eigenvalue $\Lambda^0_{0,1}(\eta)$ and the double one $\Lambda^0_{1,1}(\eta)$, the double eigenvalue $\Lambda^0_{2,1}(\eta)$ and the simple one $\Lambda^0_{0,2}(\eta)$,  the simple eigenvalue $\Lambda^0_{0,2}(\eta)$ and the double one $\Lambda^0_{3,1}(\eta)$ and finally the double eigenvalues  $\Lambda^0_{3,1}(\eta)$ and $\Lambda^0_{1,2}(\eta)$ (see Figure \ref{fig:spectrum}). The next result is a consequence of Theorem \ref{thmmain} and Corollary \ref{cor}.
\begin{cor}
   \label{cor2}
   There exists $\varepsilon_0>0$ such that for any $(n,k,m,l)\in\{(0,1,1,1),$
   \newline$(2,1,0,2), (0,2,3,1), (3,1,1,2) \}$ and for any $\varepsilon\in(0, \varepsilon_0]$, between the segments ${\cal B}_{n,k}^\varepsilon$ and ${\cal B}_{m,l}^\varepsilon$ of the spectrum of the problem \eqref{utheta}-\eqref{sbc}, there is a gap
         \begin{equation*}
         ( \overline{\Lambda}^\varepsilon_{n,k}, \underline{\Lambda}^\varepsilon_{m,l} ),
         \end{equation*}
   whose endpoints
          \begin{align*}
          \overline{\Lambda}^\varepsilon_{n,k} := \max\{\Lambda^\varepsilon_{n,k}(\eta)\hspace{0.1cm}:\hspace{0.1cm} \eta\in (-\pi, \pi]^2\}, \qquad \underline{\Lambda}^\varepsilon_{m,l} := \min\{\Lambda^\varepsilon_{m,l}(\eta)\hspace{0.1cm}:\hspace{0.1cm} \eta\in (-\pi, \pi]^2\},
          \end{align*}
   satisfy the asymptotic formulas
        \begin{equation*}
            \left|\overline{\Lambda}^\varepsilon_{n,k} - \Lambda^0_{n,k}-\varepsilon^{2m}\overline{\Lambda}^1_{n,k}  \right|\leq \tilde{C}\varepsilon^\gamma,   
             \end{equation*}
           \begin{equation*}
            \left|\underline{\Lambda}^\varepsilon_{m,l} - \Lambda^0_{m,l}- \varepsilon^{2m}\underline{\Lambda}^1_{m,l}\right|\leq \tilde{C}\varepsilon^\gamma, 
             \end{equation*}
         where $\tilde{C}=\max\{C_{n,k}, C_{m,l} \}$ and
          \begin{align*}
                   \overline{\Lambda}^1_{n,k} := \max\{\Lambda^1_{n,k}(\eta)\hspace{0.1cm}:\hspace{0.1cm} \eta\in (-\pi, \pi]^2\}, \qquad \underline{\Lambda}^1_{m,l} := \min\{\Lambda^1_{m,l}(\eta)\hspace{0.1cm}:\hspace{0.1cm} \eta\in (-\pi, \pi]^2\}.
                   \end{align*}
   \end{cor}
Now, consider $\Lambda^\varepsilon_{1,1}(\eta)$ and $\Lambda^\varepsilon_{2,1}(\eta)$. Thanks to Corollary \ref{cor}, we have that 
     \begin{equation*}
     \Lambda^\varepsilon_{1,1}(\eta)\leq \Lambda^0_{1,1}+C_{1,1}\varepsilon^\gamma,
     \end{equation*}
     \begin{equation*}
     \Lambda^\varepsilon_{2,1}(\eta)\geq \Lambda^0_{2,1}-C_{2,1}\varepsilon^\gamma.
     \end{equation*}
Since $\Lambda^0_{1,1}<\Lambda^0_{2,1}$, for small $\varepsilon$ there exists a gap between the bands $\mathcal{B}^\varepsilon_{1,1}$ and $\mathcal{B}^\varepsilon_{2,1}$. However, since the coefficients of $\varepsilon^{2m}$ vanish, one should provide more terms of the asymptotic expansion of the eigenvalues $\Lambda^\varepsilon_{1,1}(\eta)$ and $ \Lambda^\varepsilon_{2,1}(\eta)$ as well as numeric computations in order to have more information about the length of the spectral gap. Moreover, if we consider further entries in the sequence \eqref{truespectrum}, {\it i.e.}
   \begin{equation*}
   ...<\Lambda^\varepsilon_{1,2c}(\eta)\leq\Lambda^\varepsilon_{1,2s}(\eta)<\Lambda^\varepsilon_{4,1c}(\eta)\leq\Lambda^\varepsilon_{4,1s}(\eta)<\Lambda^\varepsilon_{2,2c}(\eta)\leq\Lambda^\varepsilon_{2,2s}(\eta)<\dots,
   \end{equation*}
formula \eqref{multipleeig} combined with Corollary \ref{cor} does not enable us to conclude the existence  of spectral gaps between the bands ${\cal B}^\varepsilon_{1,2}$ and ${\cal B}^\varepsilon_{4,1}$ as well as between ${\cal B}^\varepsilon_{4,1}$ and ${\cal B}^\varepsilon_{2,2}$. Indeed, we have that
    \begin{align*}
    \Lambda^\varepsilon_{1,2}(\eta)&\leq \Lambda^0_{1,2} + C_{1,2}\varepsilon^\gamma,\\
    \Lambda^0_{4,1}-C_{4,1}\varepsilon^\gamma\leq\Lambda^\varepsilon_{4,1}(\eta)&\leq \Lambda^0_{4,1} + C_{4,1}\varepsilon^\gamma,\\
    \Lambda^0_{2,2} + C_{2,2}\varepsilon^\gamma\leq\Lambda^\varepsilon_{2,2}(\eta),&
    \end{align*}
which show that a further investigation of higher order terms in the asymptotic expansion of $\Lambda^\varepsilon_{1,2}(\eta)$, $\Lambda^\varepsilon_{4,1}(\eta)$ and $\Lambda^\varepsilon_{2,2}(\eta)$ is requested in order to detect the opening of spectral gaps and it is left as an open question.
\bigskip
\noindent
\section*{Acknowledgments}
\smallskip
The authors wish to thank V. Chiad\`{o} Piat for her careful reading and relevant remarks.\\
The author S.A. Nazarov has been supported by the grant 18-01-00325 of the Russian Foundation on Basic Research. The author L. D'Elia is a member of the Gruppo Nazionale per l'Analisi Matematica, la Probabilit\`a e le loro Applicazioni (GNAMPA) of the Istituto Nazionale di Alta Matematica (INdAM).

\bibliographystyle{plain}
\bibliography{Gelfandtransform}

\begin{thebibliography}{10}

\bibitem{BK78}
T.~C. Benton and H.~D. Knoble.
\newblock Common zeros of two {B}essel functions.
\newblock {\em Math. Comp.}, 32(142):533--535, 1978.

\bibitem{BS87}
M.~Sh. Birman and M.~Z. Solomjak.
\newblock {\em Spectral theory of selfadjoint operators in {H}ilbert space}.
\newblock Mathematics and its Applications (Soviet Series). D. Reidel
  Publishing Co., Dordrecht, 1987.
\newblock Translated from the 1980 Russian original by S. Khrushch\"{e}v and V.
  Peller.

\bibitem{B58}
F.~Bowman.
\newblock {\em Introduction to {B}essel functions}.
\newblock Dover Publications Inc., New York, 1958.

\bibitem{CCNT}
A.~Cancedda, V.~Chiad\`o~Piat, S.A. Nazarov, and J.~Taskinen.
\newblock Spectral gaps for the linear water-wave problem in a channel with
  thin structures.
\newblock Submitted.

\bibitem{CDN}
V.~Chiad\`o~Piat, L.~D'Elia, and S.A. Nazarov.
\newblock The stiff {N}eumann problem: asymptotic specialty and kissing
  domains.
\newblock preprint 2019, https://arxiv.org/abs/2001.11332.

\bibitem{CPNR13}
V.~Chiad\`o~Piat, S.~A. Nazarov, and K.~Ruotsalainen.
\newblock Spectral gaps for water waves above a corrugated bottom.
\newblock {\em Proc. R. Soc. Lond. Ser. A Math. Phys. Eng. Sci.},
  469(2149):20120545, 17, 2013.

\bibitem{FT15}
F.~Ferraresso and J.~Taskinen.
\newblock Singular perturbation {D}irichlet problem in a double-periodic
  perforated plane.
\newblock {\em Ann. Univ. Ferrara Sez. VII Sci. Mat.}, 61(2):277--290, 2015.

\bibitem{FS08}
L.~Friedlander and M.~Solomyak.
\newblock On the spectrum of narrow periodic waveguides.
\newblock {\em Russ. J. Math. Phys.}, 15(2):238--242, 2008.

\bibitem{G50}
I.~M. Gelfand.
\newblock Expansion in characteristic functions of an equation with periodic
  coefficients.
\newblock {\em Doklady Akad. Nauk SSSR (N.S.)}, 73:1117--1120, 1950.

\bibitem{HL00}
R.~Hempel and K.~Lienau.
\newblock Spectral properties of periodic media in the large coupling limit.
\newblock {\em Comm. Partial Differential Equations}, 25(7-8):1445--1470, 2000.

\bibitem{HP01}
R.~Hempel and O.~Post.
\newblock Spectral gaps for periodic elliptic operators with high contrast: an
  overview.
\newblock In {\em Progress in analysis, {V}ol. {I}, {II} ({B}erlin, 2001)},
  pages 577--587. World Sci. Publ., River Edge, NJ, 2003.

\bibitem{K95}
T.~Kato.
\newblock {\em Perturbation theory for linear operators}.
\newblock Classics in Mathematics. Springer-Verlag, Berlin, 1995.
\newblock Reprint of the 1980 edition.

\bibitem{K82}
P.~A. Kuchment.
\newblock Floquet theory for partial differential equations.
\newblock {\em Uspekhi Mat. Nauk}, 37(4(226)):3--52, 240, 1982.

\bibitem{K93}
P.~A. Kuchment.
\newblock {\em Floquet theory for partial differential equations}, volume~60 of
  {\em Operator Theory: Advances and Applications}.
\newblock Birkh\"{a}user Verlag, Basel, 1993.

\bibitem{N08Neu}
S.~A. Nazarov.
\newblock Neumann problem in angular regions with periodic and parabolic
  perturbations of the boundary.
\newblock {\em Transactions of the Moscow Mathematical Society}, 69:153--208,
  2008.

\bibitem{N10-02}
S.~A. Nazarov.
\newblock Sufficient conditions for the existence of trapped modes in problems
  of the linear theory of surface waves.
\newblock {\em Zap. Nauchn. Sem. St.-Petersburg Otdel. Mat. Inst. Steklov.},
  369:pp.202--223, 2009, (English transl.: Journal of Math. Sci., 2010. V. 167,
  N 5. P. 713--725).

\bibitem{N10}
S.~A. Nazarov.
\newblock A gap in the essential spectrum of an elliptic formally selfadjoint
  system of differential equations.
\newblock {\em Differ. Uravn.}, 46(5):726--736, 2010 (English transl.:
  Differential equations. 2010. V. 46, N 5. P. 730--741).

\bibitem{N10-1}
S.~A. Nazarov.
\newblock Opening of a gap in the continuous spectrum of a periodically
  perturbed waveguide.
\newblock {\em Mat. Zametki. 2010. V. 87, N 5.}, pages 764--786, 2010(English
  transl.: Math. Notes. 2010. V. 87, N 5. P. 738--756).

\bibitem{NRT10}
S.~A. Nazarov, K.~Ruotsalainen, and J.~Taskinen.
\newblock Essential spectrum of a periodic elastic waveguide may contain
  arbitrarily many gaps.
\newblock {\em Appl. Anal.}, 89(1):109--124, 2010.

\bibitem{NRT12}
S.~A. Nazarov, K.~Ruotsalainen, and J.~Taskinen.
\newblock Spectral gaps in the {D}irichlet and {N}eumann problems on the plane
  perforated by a double-periodic family of circular holes.
\newblock volume 181, pages 164--222. 2012.
\newblock Problems in mathematical analysis. No. 63.

\bibitem{RS78}
M.~Reed and B.~Simon.
\newblock {\em Methods of modern mathematical physics. {IV}. {A}nalysis of
  operators}.
\newblock Academic Press [Harcourt Brace Jovanovich, Publishers], New
  York-London, 1978.

\bibitem{S85}
M.~M. Skriganov.
\newblock Geometric and arithmetic methods in the spectral theory of
  multidimensional periodic operators.
\newblock {\em Trudy Mat. Inst. Steklov.}, 171:122, 1985.

\bibitem{W95}
G.~N. Watson.
\newblock {\em A treatise on the theory of {B}essel functions}.
\newblock Cambridge Mathematical Library. Cambridge University Press,
  Cambridge, 1995.
\newblock Reprint of the second (1944) edition.

\bibitem{Z05}
V.~V. Zhikov.
\newblock Gaps in the spectrum of some elliptic operators in divergent form
  with periodic coefficients.
\newblock {\em Algebra i Analiz}, 16(5):34--58, 2004.

\end{thebibliography}
\end{document}